\documentclass{article}
\usepackage{amssymb}
\usepackage{amsmath}
\usepackage{amsfonts,amsthm}
\newtheorem{theorem}{Theorem}

\newtheorem{corollary}{Corollary}

\newtheorem*{definition}{Definition}

\newtheorem{lemma}{Lemma}

\newtheorem{proposition}{Proposition}
\newtheorem*{remark}{Remark}

\begin{document}
\renewcommand{\thefootnote}{}
\title{A product decomposition for the classical quasisimple groups}
\author{Nikolay Nikolov} \date{}
\maketitle
\begin{abstract} We prove that every quasisimple group of classical type is a product
of boundedly many conjugates of a quasisimple subgroup of type $A_n$.
\end{abstract} 

\section{Main Result and Notation}
\footnotetext{\emph{2000 Mathematics Subject Classification} primary 20D06, secondary 20D40}

Let $S$ be a quasisimple group of classical Lie type $\mathcal{X}$, that is one from 
$\{A_n, B_n, \\ C_n, D_n,$ $^2 A_n,$$^2 D_n\}$. Its classical definition is as a quotient of some 
group of linear transformations of a vector space 
over a finite field $F$ preserving a nondegenerate form. 
We shall not make use of this geometry but instead rely on the Lie theoretic approach. 
In other words we view $S$ as the group of fixed points of certain automorphism of an 
algebraic group defined as follows: \medskip

Let $F$ be a field (finite or infinite). In order that $S$ is quasisimple when rank of $\mathcal{X}$ is small ($\leq 2$) 
some extra conditions are usually imposed on $F$, for example that it is perfect and has enough elements. 

In case $\mathcal{X}\in \{A_n,B_n,C_n,D_n\}$ let $G(-)$ be a Chevalley group of 
Lie type $\mathcal{X}$ defined over $F$.
Then $S$ is the subgroup of $G(F)$ generated by all its unipotent elements. 

In case $\mathcal{X}$ is $A_n$ or $D_n$ 
then (depending on its isogeny type) $G(-)$ has an outer \emph{graph} automorphism $\tau$ of order 2 defined by the 
symmetry of the Dynkin diagram of $\mathcal{X}$. 
We further assume that $F$ has an automorphism $x\mapsto \bar x$ of order 2, this defines 
another involutionary automorhism $f$ of $G(F)$. The group $S$ is invariant under $\tau$ and $f$. The classical group of 
type $^2 \mathcal{X}$ is then defined as the group of fixed points of $S$ 
under the Steinberg automorphism $\sigma:=\tau f$. \medskip
 
We shall not go into further details of the definition and structure of $S$. 
The reader is assumed to be familiar with Carter's 
book \cite{C} which is our main reference. Other sources are 
\cite{GLS} or Section 2 of \cite{NS2}.
\medskip

When $\mathcal{X} \not =A_n$ we shall define a certain quasisimple subgroup $S_1$ of $S$ 
which has type $A_{n_1}$ (so it is a central quotient of  $\mathrm{SL}_{n_1+1}(F_0)$ for 
an appropriate $n_1$ and a subfield of $F$).

Our main result is
\begin{theorem}\label{main}
There is a constant $M\in \mathbb{N}$ with the following property:  
If $\mathcal{X}\not = A_n$ then 
there exist $M$ conjugates $S_i=S_1^{g_i}$ of $S_1$ such that
\[S= S_1 \cdot S_{2} \cdots S_{M}:=\{ s_1\cdots s_M \ | \ s_i \in S_i \}. \]
In fact we can take $M=200$.
\end{theorem}

This result reduces many problems about decompositions of classical groups 
to those of type $A_n$ which are usually easier. Such problems have been considered in \cite{NS2}, 
\cite{Wilson}, \cite{SW} and \cite{LS}. While the conclusion of Theorem \ref{main} is not unexpected in view 
of the large size of $S_1$ compared to $|S|$ it seems to have escaped attention so far and 
the author has not been able to find a proof of it in the literature. It is possible that 
geometric methods will yield a much better bound for $M$ but we haven't been able to 
find a simpler proof than the one given below. 

Finally we remark that a similar result holds for quasisimple groups $S$ of type $A_n$ (where it is easy): 
$S$ is a product 
of at most 5 conjugates of a subgroup of type $A_{n-1}$. Using this we can add further 
conditions on the subgroup $S_1$ in Theorem \ref{main}. See Section \ref{variations} for details. 

\section{Notation and definitions}

\begin{remark} If char $F=2$ then the simple groups of types $B_n$ and $C_n$ are isomorphic. 
In this situation we shall assume that $S$ has type 
$\mathcal{X}=B_n$. This has implications for the Chevalley commutator relations 
we use, cf. the proviso in 
the statement of Lemma \ref{comrel}. \end{remark}

In case $\mathcal{X}$ is $^2A_n$ or $^2D_n$ define $F'\subset F$ to be the subfield 
of $F$ fixed by its automorphism $x \mapsto \bar x$ of order 2. 

Let $\Pi$ be a fixed set of fundamental roots in 
the root system $\Sigma$ of $S$. We assume that $n'=|\Pi|\geq 2$. 
In case $\mathcal{X}= $ $^2A_{2n'}$ the convention in \cite{GLS} is that $\Sigma$ has type $BC_{n'}$. \medskip

The root subgroups of $S$ are denoted by $X_w$ for $w \in \Sigma$ and are usually 
one-parameter with parameter $t$ ranging over either 
\begin{itemize}
\item $F'$ (if  $\mathcal{X} \in \{^2D_n, \ ^2A_{2n'-1}\}$ and $w$ is long), or

\item $F$ (otherwise). 
\end{itemize}
The only exception to this 
is when $\mathcal{X}=$ $^2A_{2n'}$ and $w$ is a short root of $\Sigma$; then $X_w$ 
is 2-parameter, as described in \cite{GLS} Table 2.4, type IV. 
In this case we consider the center of $X_w$ as 
another (one-parameter over $F'$) root subgroup $Y_{2w} \leq X_w$ associated 
to the doubled root $2w$. (The reason for this is that we wish to define in a natural way a filtration on the 
positive unipotent group $U_+$ below). \medskip 

\begin{definition}
(a) An element $a \in X_w$ is \emph{proper} if $a=X_w(t)$ or $a=Y_{2w}(t')$, 
or (if $X_w$ is 2-parameter) $a=X_w(t,t')$ with 
$t \in F\backslash \{0\}$, respectively $t' \in F'\backslash \{0\}$.

(b) Further, when $\mathcal{X}=^2D_{n'+1}$ (and so $\Sigma$ has type $B_{n'}$), for a short root $w$ we call two
root elements $a_1,a_2\in X_w$ a \emph{proper pair} if $a=X_w(t_1), a'=X_w(t_2)$ are such that $F't_1+F't_2=F$.
\end{definition}

The fundamental roots $\Pi$ determine a partition of $\Sigma= \Sigma_+ \cup \Sigma_-$ 
into the sets of the positive and negative roots. Define
\[ U_+:= \prod_{w \in \Sigma_+}  X_w, \quad U_-:= \prod_{w \in \Sigma_-}  X_w.\]

From now on we shall assume that $\mathcal{X} \not = A_n$. \medskip

Recall the standard realization of $\Sigma$ in the $n'$-dimensional euclidean space $E$ with orthonormal basis
$e_1,e_2,...,e_{n'}$:

\[ \Sigma_+=\{\pm e_i + e_j\ | \ 1\leq i<j \leq n'\}\cup \Theta, \]
where 

$\Theta$ is empty if $\Sigma$ has type $D_{n'}$, 

$\Theta=\{e_i \ |\ 1\leq i \leq n'\}$ in type $B_{n'}$, 

$\Theta= \{2e_i \ |\ 1\leq i \leq n'\}$ in type $C_{n'}$ and

 $\Theta= \{e_i,2e_i \ |\ 1\leq i \leq n'\}$ in type $BC_{n'}$.

For this choice of $\Sigma_+$ the fundamental roots $\Pi$ are
\[\Pi= \{e_{i+1}-e_i\ | \ i=1,2,..., n'-1 \} \cup \{r_0\},\]

where if $\Sigma$ is of type $D_{n'}$ then $r_0=e_1+e_2$, in case $\Sigma$ has type $B_{n'}$ or $BC_{n'}$ 
then $r_0=e_1$ and in case $\Sigma$ has type
$C_{n'}$ we have $r_0=2e_1$.

For a root $r\in \Sigma_+$ we denote the \emph{height} of $r$ by $ht(r)$. 
This means that $r$ is a sum of $ht(r)$ fundamental roots (maybe with repetitions). 

Set $\Pi_1:= \Pi \backslash \{r_0\}$ and $\Sigma_1= \Sigma \cap \mathbb{Z}\Pi_1$.  
Then $\Sigma_1$ is a root subsystem of type $A_{n_1}$ with $n_1 =n'-1$. \medskip

Let 
\[S_1:= \langle X_v \ \mid \ v \in \pm \Pi_1 \rangle, \quad U_1:= \langle X_v \ \mid \ v \in  \Pi_1 \rangle.\]
The group $S_1 $ is a Levi factor of a parabolic subgroup of $S$, 
thus it is quasisimple of Lie type $A_{n_1}$ with $n_1 \geq 1$ and $U_1$ is its positive unipotent subgroup. 
The field of definition of $S_1$ is again $F$ with the exception of type $\mathcal{X}=$$^2D_n$ when it is $F'$. 

\section{The Proof of Theorem \ref{main}}

We use the following result by M. Liebeck and L. Pyber:
\begin{theorem}[\cite{LS}, Theorem D] If $S$ is a quasisimple group of classical Lie type then
\[ S=(U_+U_-)^6U_+.\]
\end{theorem} 

\medskip 
This reduces the problem to showing that $U_+$ is contained in a product of boundedly 
many (say $M_1$) conjugates of $U_1$. Then by symmetry the same result holds for $U_-$ and we 
can take $M=13M_1$. \bigskip

Before we proceed with the proof of Theorem \ref{main} we record a result which is a consequence of the 
Chevalley commutator relations (\cite{GLS} Theorems 1.12.1 and 2.4.5):

\begin{lemma} \label{comrel} Let $u,v$ and $u+v$ be roots of the system $\Sigma$ of classical type, such that
$u \in \Sigma_1$. 

Assume that if char $F=2$ and $S$ 
has type $\mathcal{X}=C_{n}\simeq B_n$ then $u,v$ is not a pair of short 
root summing to the long root $u+v$. 

Further, let (*) be the condition that 

\emph{ (*): The Lie type $\mathcal{X}$ of $S$ is $^2D_n$, and the roots $u+v$ and $v$ are short while $u$ is long.}

Then 
\medskip

(a) If not (*) and $a \in X_v$ is \textbf{proper} then
\[ X_{u+v}\subseteq [X_u,a]\cdot Z(u,v),\]
where $Z(u,v)$ is the product of all root subgroups $X_{w}, Y_w$ with $w=iu+jv$, $i,j \in \mathbb{N}, i+j>2$. 
(It may be that $Z=\{1\}$.)
\medskip

(b) Suppose that (*) holds.  
Then if $a_1,a_2\in X_v$ is a 
\textbf{proper pair} of root elements we have
\[ X_{u+v} \subseteq [X_u,a_1] [X_u,a_2] \cdot Z(u,v), \]
where $Z(u,v)$ is as in (a).
\end{lemma}

We now define 4 elements of $U=U_+$: \medskip

For $i=1,2,\ldots n'-1$ let $a_i$ be 
a fixed proper element of the root subgroup $X_{e_i+e_{i+1}}$. 
When $\mathcal{X}$ is different from $^2D_{n'+1}$ and $^2A_{2n'}$ 
fix a proper element $b_i$ in $X_{e_i}$ or $X_{2e_i}$ as relevant. 
If $\mathcal{X}=$ $^2D_{n'+1}$ then fix a proper pair of elements $c_{i},c'_i \in X_{e_i}$
If $\mathcal{X}=$ $^2 A_{2n'}$ then fix a proper element $d_i \in X_{e_i}$ and a proper element $d'_i\in Y_{2e_i}$. 
\begin{definition} We set

\[ w_1:=a_1a_3a_5 \cdots , \quad w_2:=a_2 a_4 a_6 \cdots, \]

If $\mathcal{X}=D_{n'}$ then set $w_3=w_4=1$.

If $\mathcal{X}=$ $^2D_{n'+1}$ define $w_3=c_1 \cdots c_{n'}$, $w_4=c'_1 \cdots c'_{n'}$.

If $\mathcal{X}=$ $^2A_{2n'}$ define $w_3= d_1 \cdots d_{n'}$, $w_4=d'_1 \cdots d'_{n'}$.

In all other cases $w_3:=b_1b_2\cdots b_{n'}$ and $w_4:=1$. 
\end{definition} \bigskip

Let $\Delta=\Sigma_+ \backslash (\Sigma_{1} \cup\{r_0,2r_0\})$ and define
\[ D=\prod_{w \in \Delta} X_w.\]

Then $D$ is a normal subgroup of $U$ and we have $U=X_{r_0}\cdot U_1\cdot D$. 

In the case when $\Sigma$ has type $B_{n'}$ we 
shall need one extra conjugate of $U_1$.

 More precisely, when $\Sigma$ has type 
$B_{n'}$ define $W=X_r$ where $r=e_1+e_2$. 
In all other cases set $W=1$. Since $e_1+e_2$ has the same length as the roots in 
$\Sigma_1$ it is clear that $W$ is contained in a conjugate of 
$U_1$, say $U_1^s$. 

We shall prove 
\begin{proposition}\label{D}
\begin{equation}\label{ind} 
D=W\prod_{j=1}^4 [U_1,w_j]. 
\end{equation}
\end{proposition}
To deal with $X_{r_0}$ we consider three cases: \medskip

\textbf{(a)} When $\Sigma$ has type $D_{n'}$ the root $r_0$ has the same length as the roots in $\Sigma_1$. 
Therefore $X_{r_0}$ is contained 
in a conjugate of $U_1$. \medskip

\textbf{(b)} Suppose $\Sigma$ is of type $B_{n'}$. We can write $r_0=e_1=u+v$ where $u=e_1-e_2\in \Sigma_-$, 
$v=e_2$. Thus $u$ is conjugate to a root in $\Sigma_{1+}$ and hence $X_u\subseteq U_1^s$ for some $s\in S$. 
Note that all roots of $\Sigma$ of the form $iu+jv$ with $i,j  \in \mathbb{N}$ except $r_0$ lie in $\Delta$. By
Lemma \ref{comrel} then
\[ X_{r_0}\subseteq [X_u,a] D \] for some proper $a\in X_v$, unless (*) holds when 
\[X_{r_0}\subseteq 
[X_u,a_1][X_u,a_2] D \]
for some proper pair $a_1,a_2 \in X_v$.   \medskip

\textbf{(c)} If $\Sigma$ has type $C_{n'}$ then $r_0=2e_1$. We set $u=e_1+e_2$, $v=e_1-e_2$ 
and the argument is as above. Note that in case that $S$ has 
type $\mathcal{X}=C_n$ the characteristic of $F$ is assumed to be different from 2 
and so the conclusion of Lemma \ref{comrel} applies. 

\textbf{(d)} Finally, if $\Sigma$ has type $BC_{n'}$ then $r_0={e_1}$ and $X_{r_0}$ is 2-parameter. 
Let $r_0=e_1\in \Pi$, 
$u=e_1-e_2 \in \Sigma_-$, 
and let $g_1$ be a proper element of the 2-parameter root subgroup $X_{e_2}$ and $g_2$ be 
a proper element of $X_{e_1+e_2}$. Then
\[ Y_{2r_0} \subset X_{r_0}\subseteq [X_u,g_1]  \cdot [X_u ,g_2] \cdot D.\]

Note that, again, $X_u$ is conjugate to a root subgroup in $U_1$.
\bigskip

In conclusion, we see that in all three cases we have 
\[U_+=X_{r_0} U_1 D = X_{r_0} U_1W\cdot \prod_{j=1}^4 [U_1,w_j] \subseteq \prod_{i=1}^6 U_1^{s_i}  \cdot 
\prod_{j=1}^4(U_1 \cdot U_1^{w_j})\]
for some appropriate choice of $s_1, \ldots s_6 \in S$.
 
Hence we can take $M_1=14$ and 
Theorem \ref{main} follows with $M=13\times 14=182$.

\subsection*{Proof of Proposition \ref{D}:} 
Let $v\in \Sigma_{1+}$ and $j\in \{1,\ldots ,4\}$ be such that $[X_v,w_j] \not = 1$. 
The element $w_j$ was defined as a product
of certain proper root elements $a_i$, $b_i$ or $c_i,c'_i$ (notation depending on $j$ and the type $\mathcal{X}$ of $S$). 
It is easy to see that 
$X_v$ commutes with all except one of the 
constituents $a_i$, (resp $b_i$ or $c_i,c'_i$) of $w_j$ ($i=1,\ldots ,n'$). 
Say this constituent is $a\in X_\alpha$ for the appropriate root 
$\alpha \in \Sigma \backslash \Sigma_1$ (or a proper pair $c,c'\in X_\alpha$ if $\mathcal{X}=$$^2D_n$).
  
Let $w=w(j,v):=v+\alpha \in \Delta$. Using Lemma \ref{comrel} it easily follows that

(A) If not (*) then 
\[ [X_v, w_j]\cdot Z = [X_v,a] \cdot Z=X_w Z,\]
where $Z=Z(w)$ is the product of root subgroups in $D$ of height $> ht(w)$. \medskip

(B) If (*) holds then 
\[ [X_v,w_3] [X_v,w_4]\cdot Z =[X_v,c] [X_v,c']\cdot Z= X_wZ.\]

Moreover, the root $w \in \Delta$ together with $w_j$ 
uniquely determines $v\in \Sigma_{1+}$.  

We set $t_j(v):=ht(w)\geq 2$ and for completeness 
declare that $t_j(v)=\infty$ if $X_w$ and $w_j$ commute. \medskip

Fixing $j$ for the moment, choose an ordering of the roots $\Sigma_{1+}$ with non-increasing $t_j$ and
write $u \in U_1$ as $u=\prod_{v \in \Sigma_1}x_v$ as a product of root elements $x_v$ in that 
(non-increasing $t_j(v)$) order. Now, using
the identity $[xy,w]=[x,w]^y \cdot [y,w]$ we obtain
\begin{equation} \label{expand} [u,w_j]=\prod_{v \in \Sigma_1}[x_v,w_j]^{h_v},
\end{equation}
where $h_v$ is an element which depends on $x_w$ with $t_j(w)\leq t_j(v)$.

For $i\geq 2$ let $D(i)$ be the product of the root subgroups $X_v \subseteq D$ with $ht(v)\geq k$. Then
$D=D(2)>D(3)>\cdots >\{1\}$ is a filtration of $D$. \medskip

We prove that the identity (\ref{ind}) holds modulo $D(k)$ for each $k\geq 2$. We use induction on $k$, starting 
with $k=2$ (when it is trivial).
Assuming that (\ref{ind}) holds modulo $D(k)$, we shall prove that it holds modulo $D(k+1)$.
 
Let $\Delta_j(k)$ be the set of roots $r\in \Sigma_1$ such that $t_j(r)=k$. 
It is clear that $\Delta_j(k) \cap \Delta_{j}(k')=
\emptyset$ if $k\not = k'$. \medskip

Let $g\in D$  be arbitrary. By the induction hypothesis
\[ g= d^{-1} y_k \prod_{j=1}^4 [g_j,w_j]\]
for some $d \in D(k)$, $y_k \in W$ and $x_j\in U_1$. We may assume that $y_k=1$ unless $k=3$ and 
$g_j\in \prod_{t_j(v)<k}X_v$. Let $x_j \in \prod_{t_j(v)=k}X_v$. Then 
\[\prod_{j=1}^4 [x_jg_j,w_j]=\prod_{j=1}^4 [x_j,w_j]^{h_j} \cdot \prod_{j=1}^4 [g_j,w_j],\]
where $h_j=g_j\prod_{l<j} [g_l,w_l]$. In turn we have 
\[ x_j=\prod_{v\in \Delta_j(k)} x_{j,v} \ \textrm{ (product in the chosen order) } , \] 
where $x_{j,v}\in X_v$ and $t_j(v)=k$. Using (\ref{expand}) we have
\[\prod_{j=1}^4 [x_jg_j,w_j] =\prod_{j=1}^4 \left( \prod_{v \in \Delta_j(k)} [x_{j,v},w_j]^{h_{j,v}} \right)
\cdot \prod_{j=1}^4 [g_j,w_j],\]
for some $h_{j,v} \in U$ depending on $h_j$ and the root elements $x_{j,w}$ 
succeeding $x_{j,v}$ in the ordering. \bigskip

Now, for $v\in \Delta_j(k)$ and $x_{j,v} \in X_v$ we have $[x_{j,v},w_j]^{h_{j,v}}\equiv [x_{j,v},w_j]$ 
mod $D(k+1)$. Therefore it is enough to prove the following:

\begin{proposition}\label{step} Let $k\geq 2$ and $d\in D(k)/D(k+1)$. There exist $x_{j,v}\in X_v$ 
for each $v\in \Delta_j(k)$ and $y\in W$ (with $y=1$ unless $k=3$ and $\Sigma$ of type $B_{n'}$) such that
\[ d \equiv y \cdot \prod_{j=1}^4 \left( \prod_{v \in \Delta_j(k)} [x_{j,v},w_j] \right) \quad \mathrm{mod} \ D(k+1).\]  
\end{proposition}

Given $v \in \Delta_j(k)$ let $r=r(j,v)$ be the unique root of height $k$ such that $[X_v,w_j]=X_r$ mod $D(k+1)$.
The group $D(k)/D(k+1)$ is abelian and product of the root subgroups $X_r$ for $r \in D$ and $ht(v)=k$. 

From Lemma \ref{comrel} and our definitions of $w_i$ it is easy to see that we only need to check 
that 

(A) If not (*) then for each $r$ in $\Delta$ of height $k$ (with the exception of $r=e_1+e_2$ in type $B_{n'}$) 
there exist at least one 
$j\in \{1,\ldots ,4\}$ and $v\in \Delta_j(k)$ such that $r=r(j,v)$.  

(B) If (*) holds then for each short root $r=e_k \in \Delta$ of height $k\geq 2$ there are long roots 
$v_3 \in \Delta_3(k), v_4 \in \Delta_4(k)$ such that $r=r(3,v_3)=r(4,v_4)$.
 \medskip

Now 
Part (B) is clear. In fact this is the reason why we introduced proper pairs and the elemnts $w_3$ and $w_4$ in 
the case when $S$ has type $^2D_{n}$.

Similarly part (A) is a matter of simple verification depending on the Lie type $\mathcal{X}$ of $S$: \medskip

\textbf{Case 1:} $\mathcal{X}=D_{n'}$.
$r=e_i+e_l$, where $i<l$ and $(i,l)\not=(1,2)$. If $i>1$ take $v=e_l-e_{i-1}$ and $j=1$ if $i$ is even, 
$j=2$ if $i$ is odd. In both cases 
we have that $X_v$ commutes with all constituents of $w_j$ except $a_{i-1}$ and thus 
$[X_v,w_j]=[X_v,a_{i-1}] =X_r$ mod $D(k+1)$. Recall that $a_{i-1}$ is a proper element 
in the root subgroup of $e_{i-1}+e_i$.
 
On the other hand, if $i=1$ then $l>2$. Take $j=1$, $v=e_l-e_2$ and the argument is as above.
\medskip

\textbf{Case 2:} $\mathcal{X}=B_{n'}$ or $^2D_{n'+1}$ so its root system $\Sigma$ has type $B_{n'}$. 
The only roots in $\Delta$ not covered in Case 1 or (B) are:

$r=e_1+e_2$: this is the exception, we have introduced $W=X_r$ precisely for this root and simply 
choose the appropriate $y\in W$.

$r=e_l$, ($n'\geq l \geq 2$) in type $\mathcal{X}=B_{n'}$: Take $v=e_l-e_{l-1}$, $j=3$ and then we have 
$[X_v,b_{l-1}]=X_r$ mod $D(k+1)$. Recall that 
$b_{l-1}$ is a proper element of the root subgroup $X_{e_{l-1}}$ and 
$w_3=b_1b_2\cdots b_{n'}$. \medskip

\textbf{Case 3:} $\mathcal{X}=C_{n'}$ or $^2A_{2n'-1}$ so $\Sigma$ has type $C_{n'}$. 
The roots $e_i+e_l$ ($(i,l)\not = (1,2)$) 
are dealt with in the same way as in Case 1. The remaining ones are:

$r=e_1+e_2$:  Take $j=3$ and $v=e_2-e_1$. Then $v$ commutes with all 
$b_i$ in $w_3=b_1\cdots b_{n'}$ except $b_1$:
 the proper element of $X_{2e_2}$. Thus $[X_v,w_j]=[X_v,b_1]=X_r$ mod $D(k+1)$.

$r=2e_l$, ($n'\geq l \geq 2$): In order to obtain long roots subgroups we need that if $\mathcal{X}=C_{n'}$ 
then the characteristic of $F$ is not 2 
(as we have assumed from the start). Take $v=e_l-e_{l-1}$ and choose $j=1,2$ such that a proper element 
$a_{l-1}$ of root subgroup 
$e_{l-1}+e_l$ is a constituent of $w_j$. (i.e. take $j=1$ if $l$ is even and $j=2$ if $l$ is odd.) \medskip

\textbf{Case 4}: $\mathcal{X}=$ $^2A_{2n'}$ so $\Sigma$ has type $BC_{n'}$. Then a combination of the 
reasoning from Cases 1, 2 and 3
gives the conclusion: 

More precisely we obtain the root subgroups $X_r$ for $r=e_i+e_l$, $(i,l)\not =(1,2)$ as in Case 1.
$X_{e_1+e_2}$ is obtained from $j=4$ and $v=e_2-e-1$. 

The root subgroups $X_{e_l}, \ (l\geq 2)$ are obtained (modulo $D(k+1) \geq Y_{2e_l}$.) 
from $[X_{e_l-e_{l-1}},w_3]$ for $j=3$ and $v=e_l-e_{l-1}$ as in Case 2.

Finally the root subgroup $Y_{2e_l}$ for $l\geq 2$ is obtained with $v=e_l-e_{l-1}$ and $j=1$ ($j=2$) if $l$ is even 
(resp. odd) just as in Case 3.   
\section{Variations} \label{variations}
First we shall consider the analogue of Theorem \ref{main} when $S$ has type $A_n$, $n\geq 3$: \medskip
 
Without loss of generality we may assume that $S=\mathrm{SL}_{n+1}(F)$, acting on $V=F^{(n)}$ with 
standard basis $(v_1,\ldots ,v_{n+1})$. For $1\leq i \leq n+1$ let $V_i$ be the subspace of $V$ 
spanned by all $v_j $ with $j\not =i$ and define 
\[ S(i)=\{ g \in S \ | \ g \cdot v_{i}=v_{i},\  g \cdot V_i=V_i \} \]
It is clear that all $S(i)$ are conjugate to each other and isomorphic to $\mathrm{SL}_{n}(F)$.
\begin{lemma} \label{SL} 
\[ S=S(3)\cdot S(2)\cdot S(1)\cdot S(2)\cdot S(3).\]
\end{lemma} 

\textbf{Proof} This is well known. For completeness we sketch one argument. 

Let $g\in S$. It is easy to see that $S(3) S(2)\cdot v_1=V\backslash \{0\}$ and 
therefore there exist $a_2 \in S_2, a_{3}\in S(3)$ such that $g\cdot v_1=a_3a_2\cdot v_1$. 
Hence the matrix $g':=a_2^{-1}a_3^{-1}g$ has the transpose of $(1,0,\cdots ,0)$ as its first column. 

By right multiplication with the elementary matrices $1+\lambda_j E_{1,j}$, with $2\leq j \leq n+1$ and appropriate 
$\lambda_j \in F$ we can make the first row of $g'$ to be $(1,0,\ldots ,0)$.  As
\[ \langle  1+\lambda_j E_{1,j}\ | \quad  2\leq j \leq n+1, \  \lambda_j \in F \rangle \subseteq S(3)S(2)\]
we conclude that there exist $b_2\in S(2), b_3\in S(3)$ such that $g'b_3b_2=g''\in S(1)$. Therefore
$g=a_3a_2g''b_2^{-1}b_3^{-1}$ as required
$\square$. \medskip

Some applications of Theorem \ref{main} may need extra conditions on the subgroup $S_1$. For 
example in \cite{NS2} it is required that $S_1$ be invariant under the group $\mathcal{D}\Phi \Gamma$ 
of diagonal-field-graph autmorphisms of $S$. The only exception to this is the case when $S$ has type 
$\mathcal{X}=D_{n'}$, when $S_1$ is not preserved by the graph automorphism $\tau$ of order 2. (Since $\tau$
does not preserve the set $\Pi_1$ of fundamental roots).
In this case we define 
\[ \overline{\Pi}:=\Pi \backslash \{r_0,r_0^\tau\}=\{e_{i+1} -e_i \ | \quad i=2,3,..,n'-1 \}\] and
$\overline{S}:=\langle X_r \  | \ r \in \pm \overline{\Pi} \rangle$. 

The group $\overline{S}$ is a Levi factor of a parabolic of $S_1$ and is conjugate to the groups $S(i)$ in 
Lemma \ref{SL} (defined for $S_1$ as an image of $S=\mathrm{SL}_{n_1+1}(F)$). Then Lemma \ref{SL}
and Theorem \ref{D} imply the following
\begin{corollary} $S$ is a product of some $5M<1000$ conjugates of $\bar{S}$.
\end{corollary}

\texttt{\bigskip}

\texttt{Nikolay Nikolov}

\texttt{New College,}

\texttt{Oxford OX1 3BN,}

\texttt{UK.}

\end{document}